\documentclass[11pt]{amsart}

\usepackage{amsmath}
\usepackage{amssymb}
\usepackage{graphicx}
\newtheorem{theorem}{Theorem}[section]
\newtheorem{corollary}[theorem]{Corollary}

\newtheorem{proposition}[theorem]{Proposition}
\theoremstyle{definition}
\newtheorem{definition}[theorem]{Definition}
\newtheorem{example}[theorem]{Example}

\usepackage{epsfig}
\usepackage{psfrag}
\usepackage{url}

\newcommand{\A}{\mathbb{A}}
\newcommand{\B}{\mathbb{B}}
\newcommand{\E}{\mathbb{E}}
\newcommand{\I}{\mathbb{I}}
\newcommand{\N}{\mathbb{N}}

\newcommand{\R}{\mathbb{R}}
\newcommand{\U}{\mathbb{U}}
\newcommand{\V}{\mathbb{V}}
\newcommand{\W}{\mathbb{W}}
\newcommand{\bbs}{\mathbb{X}}
\newcommand{\X}{\mathbb{X}}
\newcommand{\Y}{\mathbb{Y}}
\newcommand{\Z}{\mathbb{Z}}

\title[Benford's Law in dynamical systems]{A brief survey of Benford's
  Law\\ in dynamical systems}

\author[Arno Berger]{Arno Berger}
 
\address[A.\ Berger]{Department of Mathematical and Statistical
  Sciences, University of Alberta, Edmonton, {\sc Canada}}
\email{{\tt berger@ualberta.ca}}

\author[Theodore P.\ Hill]{Theodore P.\ Hill}
\address[T.P.\ Hill]{School of Mathematics, Georgia Institute of
  Technology, Atlanta, {\sc USA}}
\email{\tt theodore.hill@math.gatech.edu}



\date{\today}

\begin{document}

\begin{abstract}
This article provides a brief overview on a range of basic
dynamical systems that conform to the logarithmic distribution of
significant digits known as Benford's law. As presented here, most
theorems are special cases of known, more general results about
dynamical systems whose orbits or trajectories follow this logarithmic
law, in one way or another. These results span a wide variety of
systems: autonomous and non-autonomous; discrete- and continuous-time;
one- and multi-dimensional; deterministic and stochastic. Illustrative
examples include familiar systems such as the tent map, Newton's
root-finding algorithm, and geometric Brownian motion. The treatise is
informal, with the goal of showcasing to the specialists
the generality and universal appeal of Benford's law throughout the mathematical
field of dynamical systems. References to complete proofs are provided
for each known result, while one new theorem is presented in
some detail.
\end{abstract}

\maketitle

%%%%%%%%%%%%%%%%%%%%%%%%%%%%%%%%%%%%%%%%%%%%%%%%%%%%%%%%%%%%%%%%%%%%%%%%%%%%%%%%

\section{Introduction}

The purpose of this article, written in honor of the authors' friend and
colleague Professor Leonid Bunimovich on the occasion of his 75th
birthday, is to present an informal overview of the remarkable
ubiquity of the statistical phenomenon called {\bf Benford's law}
(henceforth abbreviated as {\bf BL}) in the mathematical field of dynamical systems. 
This topic is especially fitting to honor Leonid, since it was in
joint work with him \cite{BBH} that the authors first explored some of
its many intriguing aspects, in response to the question ``Do
dynamical systems follow Benford's law?" raised by Tolle, Budzien, and
LaViolette \cite{TBL00} after they had found experimental evidence of
BL in several familiar physical dynamical systems; see also \cite
{Bia15, SamTJ10, SCD01}. 

The tone of this survey is intentionally informal with the goal of
introducing researchers in dynamical systems to a variety of
mathematical examples of BL. As such, most of the theorems presented
here are not new, nor are they stated in the most general
form possible. One new theorem (Theorem \ref{thm19a} below) is
included, however, and references to the proofs of all other results
are given. As the reader is going to see, the results 
primarily take the form of sufficient conditions for a variety of dynamical
systems to conform to BL, in one way or another. 

The article is organized as follows: Section \ref{sec2} provides a
brief introduction to BL. Sections \ref{sec3} to \ref{sec5} present a
selection of typical results regarding BL for various basic classes of
dynamical systems. Examples are given to illustrate each theorem,
though the presentation of these examples is brief, with
most of the (usually elementary) details being left to the interested reader. In
keeping with the informal tone, the selection of material in this
survey is somewhat narrow; for a more comprehensive overview the
reader may want to turn, e.g., to \cite{BerAH15}. Finally, Section
\ref{sec6} establishes conformance to BL for the solutions of a class
of highly non-linear two-step recursions, utilizing dynamical systems
ideas and techniques. 

Throughout, standard notation and terminology is used. Specifically,
$\N$, $\Z$, $\R_+$, and $\R$ denote the sets of all positive integers,
integers, positive reals, and reals, respectively, equipped with their usual
arithmetic, order, and topology. Half-open intervals in $\R$ are
written as $[a,b[\: =\! \{x\in \R : a\le x < b\}$, and similarly for open or
closed intervals. Also, $\lfloor x \rfloor = \max \{k\in \Z : k \le
x\}$ for every $x\in \R$. The cardinality of any finite set $A$ is $\#
A$. As usual, a function $T:\R \to \R$ being $C^m$, with $m\in \N$,
means that $T$ is $m$ times continuously differentiable.

%%%%%%%%%%%%%%%%%%%%%%%%%%%%%%%%%%%%%%%%%%%%%%%%%%%%%%%%%%%%%%%%%%%%%%%%%%%%%%%%

\section{Basic definitions}\label{sec2}

For the purpose of this article, by a {\em dynamical system\/} is
simply meant any mathematical model of how a particular process
unfolds over time. Given is a fixed set $\bbs$ (``phase space''),
together with a rule $T$ that specifies how the process unfolds in
$\bbs$. The phase space may consist of real numbers, vectors,
or random variables, etc., and the changes may occur in discrete or
continuous time. Given $x_0\in \bbs$ and a rule $T$, the orbit
$O_{T}(x_0)$ is the trajectory of the process as it moves through the
state space $\bbs$, starting at $x_0$ and proceeding according to the
rule specified by $T$. 

In this survey, the emphasis is on {\em decimal\/} representations of
numbers, the classical setting of BL, so here and throughout $\log t $
denotes the base-$10$ logarithm of $t\in\R_+ $, and all digits are
{\em decimal\/} digits. For other bases such as binary or hexadecimal,
analogous results hold with very little change, simply by replacing
$\log$ with $\log_b$ for the appropriate base $b>1$; the interested
reader is referred to \cite[p.\ 9]{BerAH15} for details. 

The basic notion underlying BL concerns the {\em first\/} (or {\em leading})
{\em significant digit\/} and, more generally, the {\em significand\/} of
a real number. Formally, for every $x$ in $\R_+$, the ({\bf decimal}) {\bf significand\/} of $x$, denoted
$S(x)$, is given by $S(x) = t$, where $t$ is the unique number in $[1,
10[$ with $x = 10^k t$ for some (necessarily unique) $k \in \Z$. For
negative $x$, let $S(x) = S(-x)$, and also $S(0)=0$ for
convenience. The integer $D_1(x) = \lfloor S(x) \rfloor$ is the {\bf
  first significant} ({\bf decimal}) {\bf digit} of $x$, so
$D_1(x) \in \{1,2,\ldots , 9\}$ for every $x$ in $\R \setminus \{0\}$.  
For example, $S(2025) = 2.025 = S(0.02025) =S(-20.25) $ and
$D_1(2025) = 2 = D_1(0.02025) = D_1(-20.25)$.

\begin{figure}[ht]
\psfrag{taa}[]{$D_1$}
\psfrag{td1}[]{$1$}
\psfrag{td2}[]{$2$}
\psfrag{td3}[]{$3$}
\psfrag{td4}[]{$4$}
\psfrag{td5}[]{$5$}
\psfrag{td6}[]{$6$}
\psfrag{td7}[]{$7$}
\psfrag{td8}[]{$8$}
\psfrag{td9}[]{$9$}

\psfrag{tnc1}[r]{\small $30.10$}
\psfrag{tnc2}[r]{\small $17.61$}
\psfrag{tnc3}[r]{\small $12.49$}
\psfrag{tnc4}[r]{\small $9.70$}
\psfrag{tnc5}[r]{\small $7.91$}
\psfrag{tnc6}[r]{\small $6.70$}
\psfrag{tnc7}[r]{\small $5.79$}
\psfrag{tnc8}[r]{\small $5.12$}
\psfrag{tnc9}[r]{\small $4.58$}

\psfrag{tnd1}[r]{\small $30.11$}
\psfrag{tnd2}[r]{\small $17.62$}
\psfrag{tnd3}[r]{\small $12.50$}
\psfrag{tnd4}[r]{\small $9.68$}
\psfrag{tnd5}[r]{\small $7.92$}
\psfrag{tnd6}[r]{\small $6.68$}
\psfrag{tnd7}[r]{\small $5.80$}
\psfrag{tnd8}[r]{\small $5.13$}
\psfrag{tnd9}[r]{\small $4.56$}

\psfrag{tne1}[r]{\small $29.56$}
\psfrag{tne2}[r]{\small $17.89$}
\psfrag{tne3}[r]{\small $12.76$}
\psfrag{tne4}[r]{\small $9.63$}
\psfrag{tne5}[r]{\small $7.94$}
\psfrag{tne6}[r]{\small $7.15$}
\psfrag{tne7}[r]{\small $5.71$}
\psfrag{tne8}[r]{\small $5.10$}
\psfrag{tne9}[r]{\small $4.26$}

\psfrag{tng1}[r]{\small $30.10$}
\psfrag{tng2}[r]{\small $17.60$}
\psfrag{tng3}[r]{\small $12.49$}
\psfrag{tng4}[r]{\small $9.69$}
\psfrag{tng5}[r]{\small $7.91$}
\psfrag{tng6}[r]{\small $6.69$}
\psfrag{tng7}[r]{\small $5.79$}
\psfrag{tng8}[r]{\small $5.11$}
\psfrag{tng9}[r]{\small $4.57$}

\psfrag{tcc1}[]{$(2^n)$}
\psfrag{tdd1}[]{$(F_n)$}
\psfrag{tee1}[]{$(n!)$}
\psfrag{tgg1}[]{\small  exact BL} 

\begin{center}
\includegraphics{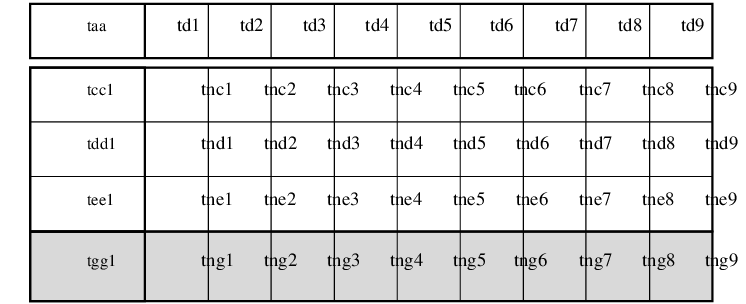}
%
%  magnification 1.00
%
\caption{Proportions (in percent) of first
  significant digits $D_1=1,2,\ldots, 9$ for the first $N= 10^4$ entries
  of $(x_n)$, for the {\em first-order autonomous\/} system $x_n = 2
  x_{n-1}$ and $x_0=1$ (powers of $2$); the {\em second-order linear
    system\/} $x_n = x_{n-1} + x_{n-2}$ and $x_1=x_2=1$ (Fibonacci
  numbers); and the {\em non-autonomous\/} system $x_n = n x_{n-1}$
  and $x_0=1$ (factorials). Each proportion matches quite closely the
  corresponding precise proportion in BL (bottom row), in the case of
  $(2^n)$ and $(F_n)$ quite strikingly so; see also Example
  \ref{t14e1a} below. This match is indicative of the fact
  that all three sequences indeed are Benford, according to Definition
  \ref{def1}.}\label{figure1}
\end{center}
\end{figure}

Informally, a real-valued orbit, trajectory, or dataset is {\em Benford\/} if
its significands are logarithmically distributed, that is, if for every
$t \in [1, 10[$ the (limiting) proportion of significands
with value $\leq t$ is exactly $\log t$.  For instance, if a sequence
of numbers is Benford (precise definition below), then exactly $\log 2
\cong 30.10\%$ of these numbers have first significant digit $D_1=1$,
exactly $\log \frac32  \cong 17.60\%$ have first significant digit $D_1=2$,
etc. Figure \ref{figure1} shows the statistics of $D_1$ for three
familiar integer sequences. 
As suggested by this figure, it is straightforward to perform a
preliminary check on whether or not a dataset is Benford, simply by
looking at the frequencies of the first significant digits.  For this
reason, BL or, more precisely, the {\em first-digit law\/} derived from
it, is widely used to check for indications of
anomalies in data; see, e.g.,  Chapter 10 in \cite{BerAH15}, or search \cite{BerAHR09} for ``fraud".  

The following definition formalizes conformance to BL for various
mathematical objects that may arise from dynamical systems.

\begin{definition}\label{def1} 
\begin{enumerate}
\item
A (real) {\em sequence\/} $(x_n)= (x_1, x_2, x_3, \ldots)$ is {\bf Benford} if
$$
\lim\nolimits_{N \to \infty} \frac{\# \{ 1 \le  n \le N : S(x_n) \le t \}}{N} = \log t 
\qquad \forall  t \in [1,10[\, ;
$$ 
\item
a (Borel) {\em function\/} $f : [0, \infty[\:  \rightarrow \R$ is {\bf
  Benford\/} if, with $\lambda$ denoting Lebesgue measure,
$$
 \lim\nolimits_{N \to \infty} \frac{\lambda \bigl( \bigl\{ x \in [0,
   N[\:  : S\bigl( f(x) \bigr)
     \leq t \bigr\} \bigr) }{N} = \log t \qquad \forall   t \in [1,10[\, ; 
$$
\item
a (real-valued) {\em random variable\/} $X$ is {\bf Benford\/} if 
$$
P\{ S( X ) \leq t\}  = \log t \qquad \forall t \in [1, 10[\, ; 
$$
\item
a {\em sequence of random variables\/} $(X_n)= (X_1, X_2, X_3, \ldots)$ {\bf converges in
  distribution to BL\/} if
$$
\lim\nolimits_{n\to \infty} P\{ S(X_n) \leq t\} = \log t \qquad
\forall t  \in [1,10[\, ,
$$
and is {\bf Benford with probability one\/} if 
$$
P\{ (X_1, X_2, X_3, \ldots)  \: \mbox{\rm is Benford} \, \} = 1; 
$$
\item 
a {\em stochastic process\/} $(X_t)_{t\ge 0}$ {\bf converges weakly to
  BL} if $(X_{t_n})$ converges
in distribution to BL for every increasing sequence $(t_n)$ with $t_n
\to \infty$, and is {\bf Benford with probability one} if the function
(or {\em path\/}) $t\mapsto X_t$ is Benford with probability one.
\end{enumerate}
\end{definition}

As an immediate consequence of Definition \ref{def1}, if a sequence
$(x_n)$ is Benford then
\begin{equation}\label{eqNxN1}
\lim\nolimits_{N\to \infty} \frac{\# \{1\le n \le N : D_1(x_n) =
  j\}}{N} = \log \frac{j+1}{j}\quad  \: \: \: \forall j \in
\{1,2,\ldots , 9\} \, ,
\end{equation}
and similarly for a Benford function, random variable etc. Though BL
is most often encountered in this first-digit form, notice that
(\ref{eqNxN1}) is weaker than the Benford property; for instance, the
sequence $(x_n) = \bigl( D_1(2^n)\bigr)$ satisfies (\ref{eqNxN1}) but
clearly is not Benford.

Many intriguing and not immediately obvious characterizations for the
Benford properties of Definition \ref{def1} can be found in the
literature; see, e.g., Chapters 4 and 5 of \cite{BerAH15}. For later
use, only one such characterization pertaining to property (i) is
stated here; throughout, let $\log 0 = 0$ for convenience, and
abbreviate {\em uniformly distributed modulo one\/} as {\em u.d.\
  mod\/} 1.

\begin{theorem}\label{thm1a}
  For every sequence $(x_n)$ in $\R$ the following are equivalent:
  \begin{enumerate}
  \item $(x_n)$ is Benford;
  \item $(ax_n)$ is Benford for every $a\in \R \setminus \{0\}$;
  \item $(\log |x_n|)$ is u.d.\ mod $1$.
  \end{enumerate}
\end{theorem}

\begin{proof}
See \cite[Thms.\ 4.2 and 4.4]{BerAH15}.
\end{proof}

%%%%%%%%%%%%%%%%%%%%%%%%%%%%%%%%%%%%%%%%%%%%%%%%%%%%%%%%%%%%%%%%%%%%%%%%%%%%%%%%

\section{One-dimensional discrete-time systems}\label{sec3}

Throughout, given any set $\bbs$, any map $T:\bbs \to \bbs$ or sequence of
maps $T_n:\bbs \to \bbs$, and any point $x_0 \in \bbs$, denote by
$O_T(x_0)$ the unique sequence in $\bbs$ with $x_n = T(x_{n-1})$ or
$x_n = T_n (x_{n-1})$ for all $n\in \N$. Refer to $O_T(x_0)$ as the
{\bf orbit} of $x_0$ under $T$ or $(T_n)$. A key object of interest
throughout this article is what may be referred to as the {\em Benford
  set\/} $\B$: If $\X \subset \R$ then simply
$$
\B  = \bigl\{  x_0 \in \X : O_T(x_0) \: \mbox{\rm is Benford}\bigr\}
\, .
$$
If $\B$ constitutes a (measure-theoretically or topologically) large part, or even all of $\X$ then BL can
rightly be viewed as a prominent feature of the dynamics of $T$ or $(T_n)$.

\subsection{Autonomous systems}

\begin{theorem} \label{t1}
Let $T:\R_+ \rightarrow \R_+$ satisfy $T(x) = ax + g(x)$ for some $a >
1$ and $g=g(x)\ge 0$.
If $\log a$ is irrational and $g(x)=o(x)$ as $x\to \infty$ then $\B = \R_+$;
if $\log a$ is rational and $g(x) = o(x/\log x)$ as $x\to
  \infty$ then $\B = \varnothing$.
\end{theorem}
\begin{proof}
See \cite[Thm.\ 6.13]{BerAH15}.
\end{proof}

\begin{example}\label{t1e1}
The orbit of $x_0$ under $T(x) = ax$ is $(ax_0, a^2 x_0, a^3 x_0,
\ldots)$. By Theorem \ref{t1}, the sequences $(2^n)$ and $(5^n)$ are
Benford; similarly $(5\cdot 2^n )$ and $(2\cdot 5^n)$ are Benford. Clearly $(10^n)$ is not Benford.
\end{example}

\begin{example}\label{t1e2}
By Theorem \ref{t1}, the orbit of $x_0$ under $T(x) = 2x + e^{-x}$ is
Benford for every $x_0 \in \R_+$, and indeed for every $x_0\in \R$ ---
notwithstanding the fact that no simple explicit formula for $T^n$ is known to the authors.
\end{example}

\begin{theorem} \label{t2}
Let $T:\R \rightarrow \R$  be $C^2$ with $T(0) = 0 $ and $0 < |T'(0)|
< 1$, and $c\in \R_+$ sufficiently small. If $\log |T'(0)|$ is
irrational then $\B \supset [-c,c] \setminus \{0\} $; if $\log |T'(0)|$ is
rational then $\B \cap [-c,c ]  =\varnothing$.
\end{theorem}
\begin{proof}
See \cite[Cor.\ 6.18]{BerAH15}.
\end{proof}

\begin{example}\label{t2e1}
Given any $a\in \R$, the map $T(x)= x + ae^{-x} - a$ is smooth with
$T(0) = 0$ and $T^{'}(0) = 1-a$. By Theorem \ref{t2}, the orbit
$O_{T}(x_0)$ is Benford for all $x_0$ sufficiently close to $0$ if $a=
0.1$, but for no such $x_0$ if $a = 0.9$. In fact, $O_T(x_0)$ is
Benford unless $T(x_0)=0$ if $a=0.1$, but is not Benford for any
$x_0\in \R$ if $a=0.9$.
\end{example}

For the statement of the next theorem, as well as several results in a
similar spirit later on, say that $\B$ has {\bf full measure} in a set
$\Y \subset \X$ if $\Y \setminus \B$ is a nullset. Thus,
if $\B$ has full measure in $\Y$ then (Lebesgue) almost every element of $\Y$
also belongs to $\B$. Say that $\B$ is {\bf meagre} (or {\bf of first category}) in $\Y$ if
$ \B\cap \Y$ is a union of countably many nowhere dense sets in $\Y$. Note
that if $\B$ is meagre in $\Y$, and if $\Y$ is a Baire space (with the
topology inherited from $\X$) with no isolated points, e.g., if $\Y$ is
any non-degenerate interval, then $\Y \setminus \B$ is dense in $\Y$
and uncountable; in particular, $\Y\setminus \B \ne \varnothing$ in
this case. As is well known, it is not uncommon in dynamics for sets
of interest (or their complements) to have full measure and yet be
meagre. The following simple result indicates how such a situation may occur
naturally in any study of BL based on Theorem \ref{thm1a}; see,
e.g., \cite[Lem.\ 3.2]{B-Kloeden60}.

\begin{proposition}\label{prop25a}
  Given any number $b >1$ and sequence $(c_n)$ in $\R$, the set
  $$
  \bigl\{
y\in \R : (b^n y + c_n) \: \mbox{\rm is u.d.\ mod 1} 
  \bigr\}
  $$
has full measure but is meagre in $\R$.
\end{proposition}

\begin{theorem} \label{t3}
Let $T:\R_+\rightarrow\R_+$ be $C^1$ with $T(x) = ax^b + g(x)$ for
some $a\in \R_+, b>1$. If $g'(x) =
o(x^{b-1}/\log x)$ as $x\to \infty$ then there exists $c\in \R_+$ so that
$\B$ has full measure but is meagre in $[c,\infty[$.
\end{theorem} 
\begin{proof} 
See \cite[Thm.\ 6.23]{BerAH15}.
\end{proof}

In the context of Theorem \ref{t3}, note that while the set $\B
\cap[c,\infty[$ is large measure-theoretically, the set
$[c,\infty[\setminus \B$ is large topologically as it is {\em
  residual}, i.e., it contains an intersection of countably many open
dense sets, and hence is dense in $[c,\infty[$ and uncountable.
Thus the Benford set and its complement
are intertwined in a very intricate way indeed. Such a state of
affairs is going to be observed repeatedly throughout the remainder of this article.

\begin{example}\label{t3e1}
The orbit of $x_0$ under $T(x) = x^2 + 1$ is Benford for almost all
$x_0 \in \R$, but the set $\{ x_0 \in \R : O_T(x_0) \: \mbox{\rm is not Benford}
\}$ is residual in $\R$, hence in particular dense and uncountable. It is unknown whether or not the integer
sequence $O_T(1)$, or in fact $O_T(k)$ for any $k\in \Z$, is Benford; see \cite[Exp.\ 6.25]{BerAH15}.
\end{example}

\begin{theorem} \label{t4}
Let $T:\R \rightarrow \R$ be real-analytic with $T(0) = T^{'}(0) =
0$. If $T$ is not constant then there exists $c\in \R_+$ so that $\B$
has full measure but is meagre in $[-c,c]$.
\end{theorem} 
\begin{proof}
See \cite[Cor.\ 6.28]{BerAH15}.
\end{proof}

\begin{example}\label{t4e1}
The orbit of $x_0$ under $T(x) = x-1 + e^{-x}$ is Benford for almost
all, but not all, $x_0\in \R$.
\end{example}

As the reader will have noticed, the previous theorems address Benford
behavior of orbits generated by maps that grow (at most) polynomially as $x\to
\infty$ or $x\to 0$.  The next result presents a simple condition guaranteeing
Benford behavior in orbits under exponential functions.

\begin{theorem} \label{t5}
Let $T:\R_+ \rightarrow \R_+$ be such that $x\mapsto \log T(10^x)$ is
convex. If $T(10^x)> 10^x T(1)$ for all $x\in \R_+$ then $\B$ has full
measure but is meagre in $\R_+$. 
\end{theorem} 
\begin{proof}
See \cite[Prop.\ 6.31]{BerAH15}.
\end{proof}

\begin{example}\label{t5e1}
The orbit $O_{T}(x_0)$ of $x_0$ under $T(x) = e^x$ is Benford for
almost all, but not all, $x_0 \in \R$.
\end{example}

While the results in this section so far pertain to systems that are
rather unexciting dynamically, the following is about a classical
``chaotic'' system.

\begin{theorem} \label{t9e21}
  Let $T:\R \rightarrow \R$ be the tent map $T(x) = 1-| 2x - 1|$. Then
$\B\supset \R \setminus [0,1]$, and $\B\cap [0,1]$ is a dense, uncountable nullset
in $[0,1]$.
\end{theorem}

\begin{proof}
See \cite[Exp.\ 6.52]{BerAH15}.
\end{proof}

One of the most widely used numerical root-finding algorithms is
{\em Newton's method}, and as the next theorem shows, the differences
between successive Newton approximations conform to BL, and so do the
differences between these approximations and the (usually) unknown root.

\begin{theorem} \label{t6}
  Let $\I \subset \R$ be a non-empty open interval, and assume that
  $f:\I\rightarrow \R$ is real-analytic but not linear, with $f(x^*) =
  0$ for some $x^* \in \I$. Let $T:\I \rightarrow \R$ be the Newton
  root-finding map associated with $f$, that is, $T(x) = x- f(x)/f'(x)$ if
  $f'(x)\ne 0$ and $T(x) = x$ if $f'(x) = 0$.
\begin{enumerate}
\item If $x^*$ is a simple root of $f$ then the sequences
  $\bigl(T^{n+1}(x_0) - T^{n}(x_0)\bigr)$  and $(T^{n}(x_0)
  - x^*)$ both are Benford for almost all, but not all, $x_0$ in a neighbourhood of $x^*$.
\item If $x^*$ is a root of $f$ with multiplicity at least $2$ then the sequences
  $\bigl( T^{n+1}(x_0) - T^{n}(x_0)\bigr)$ and $(T^{n}(x_0) - x^*)$ both are Benford
  for all $x_0 \neq x^*$ in a neighbourhood of $x^*$.
\end{enumerate} 
\end{theorem}

\begin{proof}
See \cite[Thm.\ 6.35]{BerAH15}.
\end{proof}

\begin{example}\label{t6e1}
The Newton approximations associated with $f(x) = e^x - 2$ and $g(x) = (e^x
-2)^3$ both are Benford in the sense of Theorem \ref {t6}.
\end{example}

\subsection{Non-autonomous systems}

\begin{theorem} \label{t7}
For every $n \in \N$, let $T_n:\R_+ \rightarrow \R_+$ be given by
$T_{n}(x) = a_{n}x$ for some $a_n \in \R_+$. If
$O_T(1) = (a_1, a_1a_2, a_1a_2a_3, \ldots)$ is Benford then $\B = \R_+$; otherwise $\B
= \varnothing$.
\end{theorem}

\begin{proof}
See \cite[Thm.\ 6.40]{BerAH15}.
\end{proof}

\begin{example}\label{t7e1}
Let $T_n(x) = (2 + 1/n) x$ for each $n \in \N$. By Theorem \ref{t7},
$O_T(x_0)$ is Benford for all $x_0 \in \R_+$ since, as can be shown,
the sequence $O_T(1)= (3, 3 \cdot 5/2!,  3\cdot 5\cdot 7/ 3!, \ldots)$ is Benford.
\end{example}

\begin{theorem} \label{t8}
For every $n \in \N$, let $T_n:\R_+ \rightarrow \R_+$ be given by
$T_{n}(x) = a_{n}x^{b_n}$ for some $a_n \in \R_+, b_n \in \R \setminus
\{0\}$. If $\liminf_{n \to \infty} | b_n| > 1$ then $\B$ has full
measure but is meagre in $\R_+$.
\end{theorem}

\begin{proof}
See \cite[Thm.\ 6.46]{BerAH15}.
\end{proof}

\begin{example}\label{t8e1}
Let $T_n(x) = 2^nx^2$ for every $n \in \N$.  By Theorem \ref{t8},
$O_{T}(x_0) $ is Benford for almost all, but not all, $x_0 \in \R$. 
\end{example}

\begin{theorem} \label{t9}
For every $n \in \N$, let $T_n:\R_+ \rightarrow \R_+$ be such that
$x\mapsto \log T_n (10^x)$ is convex on $\R_+$, and $T_n(10^x) \ge
10^{b_n x} T_n(1)$ for some $b_n\in \R_+$ and all $n\in \N$, $x\in
\R_+$. If $\liminf_{n\to \infty} {b_n} > 1$  then $\B$ has full
measure but is meagre in $[c,\infty[$ for some $c\in \R_+$.
\end{theorem}

\begin{proof}
See \cite[Thm.\ 3.3]{B-Kloeden60} and \cite[Thm.\ 6.49]{BerAH15}.
\end{proof}

\begin{example}\label{t9e1}
Let $T_n(x) = x^{2n} + 1$ for every $n \in \N$.  By Theorem \ref{t9},
$O_{T}(x_0)$ is Benford for almost all, but not all, $x_0\in \R$.
\end{example}

%%%%%%%%%%%%%%%%%%%%%%%%%%%%%%%%%%%%%%%%%%%%%%%%%%%%%%%%%%%%%%%%%%%%%%%%%%%%%%%%

\section{One-dimensional continuous-time dynamical systems}\label{sec4}

This short section complements the results of the previous section by
considering two continuous-time scenarios. Given any continuously
differentiable function $T:\R\to \R$ and any point $x_0 \in \R$, assume
that the initial value problem 
$$
\dot{x} := \frac{{\rm d}x}{{\rm d}t} = T(x), \quad x(0) = x_0 ,
$$
has a unique solution $x = x(t)$ for all $t \geq 0$; to emphasize the
analogy to the previous section, denote this solution by $O_T(x_0)$.

\begin{theorem} \label{t15}
Let $T:\R \rightarrow \R$ be $C^1$ with $T(x) > 0$ for all
$x\in \R_+$. If $T(x)/x$ either converges to a (finite) positive limit
as $x \to \infty$ or is non-decreasing on $]c, \infty[$ for some
$c\ge 0$, then $\B \supset \R_+$.
\end{theorem}

\begin{proof}
See \cite[Thm.\ 6.58]{BerAH15}.
\end{proof}

\begin{example} \label{t15e1}
Clearly the solution $x(t) = x_0 e^t$ of the initial value problem
$\dot x = x$, $x(0)=x_0$ is Benford for all $x_0 \ne 0$. Similarly, by
Theorem \ref{t15}, the solution of $$\dot x = \sqrt{x^2 + 1}, \quad  x(0) = x_0
$$ 
is Benford for every $x_0 \in \R$. 
\end{example}
 
\begin{theorem} \label{t16}
Let $T:\R \rightarrow \R$ be $C^1$ with $T(0) = 0$. If
$T^{'}(0) < 0$ then $\B \supset [-c,c]\setminus \{0\} $ for some $c\in
\R_+$.
\end{theorem}

\begin{proof}
See \cite[Thm.\ 6.60]{BerAH15}.
\end{proof}

\begin{example} \label{t16e1}

For every $x_0 \neq 0$, the solution of 
$$\dot x  = - \frac{x}{x^2 +1}, \quad x(0) = x_0
$$ 
eventually gets close to 0, and hence is Benford by Theorem \ref{t16}.
\end{example}

%%%%%%%%%%%%%%%%%%%%%%%%%%%%%%%%%%%%%%%%%%%%%%%%%%%%%%%%%%%%%%%%%%%%%%%%%%%%%%%%

\section{Multi-dimensional and stochastic systems}\label{sec5}

This section presents three scenarios where the phase space $\bbs$ is
not one-dimensional, and may in fact be infinite-dimensional. Though
they are stated here in their most elementary and concrete forms, the results in
all three scenarios can also be understood more abstractly as pertaining
to the dynamics of $T:\bbs \to \bbs$ given by $T(x) = x \bullet x_0$
for all $x\in \bbs$, where $x_0 \in \bbs$ is fixed and $\bullet$
denotes an appropriate binary operation on $\bbs$, such as, e.g., the
multiplication of square matrices or the convolution of distribution
functions. As the reader will no doubt appreciate, these results
typically require for their proof more sophisticated tools than the
results in previous sections; in keeping with the informal nature of this
survey, however, proofs are again provided solely via precise references.

For a first scenario, fix $d\in \N$ and let $\bbs = \R^{d\times d}$,
the set (in fact, $\R$-algebra) of all real $d\times
d$-matrices. Given $A\in \bbs$, denote by $\rho (A)$ and $|A|$ the
spectral radius and the (spectral or Euclidean) norm of $A$,
respectively, and by $[A]_{k \ell}$ the entry of $A$ in the $k$th row and
$\ell$th column, with $k,\ell \in \{1,2, \ldots , d\}$; as ususal, $A>0$
means that $[A]_{k \ell}>0$ for all $k,\ell$.

\begin{theorem}\label{thm14}
Let $A\in \bbs$, and assume that $A^m>0$ for some $m\in \N$. Then the
following are equivalent:
\begin{enumerate}
\item $\log \rho(A)$ is irrational;
\item the sequence $(|A^n|)$ is Benford;
\item for every $k,\ell \in \{1,2,\ldots , d\}$ the sequence
  $([A^n]_{k \ell})$ is Benford.
\end{enumerate}
\end{theorem}

\begin{proof}
See \cite[Thm.\ 3.2]{BE} and \cite[Thm.\ 7.11]{BerAH15}.
\end{proof}

\begin{example} \label{t14e1}
Let $A = \left[
\begin{array}{cc} 1 & 1 \\ 1 & 0 \end{array}
    \right]$. Since $A^2 =  \left[
\begin{array}{cc} 2 & 1 \\ 1 & 1 \end{array}
    \right]>0$ and $\log \rho (A) = \log \frac12 (1+\sqrt{5})$ is
    irrational, the sequence $(|A^n|)$ is Benford, and so is
    $([A^n]_{k \ell})$ for every $k,\ell \in \{1,2\}$.
\end{example}

A simple corollary of Theorem \ref{thm14} pertains to linear
recursions.

\begin{theorem} \label{t14}
Let $d \in \N$, and $a_1, a_2, \ldots, a_d , x_1,x_2, \ldots, x_d \in
\R_+$. Assume that $(x_n)$ is the (unique) solution of
$$
x_n = a_1 x_{n-1} + a_2 x_{n-2} + \ldots + a_d x_{n-d} \qquad \forall
n \ge d+1  .
$$
Then $(x_n)$ is Benford if and only if $\log \zeta$ is irrational,
where $z=\zeta$ is the (unique) root of $z^d= a_1 z^{d-1} + a_2 z^{d-2} +
\ldots + a_d$ with the largest real part. 
\end{theorem}

\begin{proof}
See \cite[Thm.\ 7.39]{BerAH15}.
\end{proof}

\begin{figure}[ht]
\psfrag{taa}[]{$D_1$}
\psfrag{td1}[]{$1$}
\psfrag{td2}[]{$2$}
\psfrag{td3}[]{$3$}
\psfrag{td4}[]{$4$}
\psfrag{td5}[]{$5$}
\psfrag{td6}[]{$6$}
\psfrag{td7}[]{$7$}
\psfrag{td8}[]{$8$}
\psfrag{td9}[]{$9$}

\psfrag{tnc1}[]{\small $30$}
\psfrag{tnc2}[]{\small $18$}
\psfrag{tnc3}[]{\small $13$}
\psfrag{tnc4}[]{\small $9$}
\psfrag{tnc5}[]{\small $8$}
\psfrag{tnc6}[]{\small $6$}
\psfrag{tnc7}[]{\small $5$}
\psfrag{tnc8}[]{\small $7$}
\psfrag{tnc9}[]{\small $4$}

\psfrag{tnc1a}[]{\small ${\it 30}$}
\psfrag{tnc2a}[]{\small ${\it 18}$}
\psfrag{tnc3a}[]{\small ${\it 12}$}
\psfrag{tnc4a}[]{\small ${\it 10}$}
\psfrag{tnc5a}[]{\small ${\it 8}$}
\psfrag{tnc6a}[]{\small ${\it 7}$}
\psfrag{tnc7a}[]{\small ${\it 6}$}
\psfrag{tnc8a}[]{\small ${\it 5}$}
\psfrag{tnc9a}[]{\small ${\it 4}$}

\psfrag{tnd1}[]{\small $301$}
\psfrag{tnd2}[]{\small $177$}
\psfrag{tnd3}[]{\small $125$}
\psfrag{tnd4}[]{\small $96$}
\psfrag{tnd5}[]{\small $80$}
\psfrag{tnd6}[]{\small $67$}
\psfrag{tnd7}[]{\small $56$}
\psfrag{tnd8}[]{\small $53$}
\psfrag{tnd9}[]{\small $45$}

\psfrag{tnd1a}[]{\small ${\it 301}$}
\psfrag{tnd2a}[]{\small ${\it 176}$}
\psfrag{tnd3a}[]{\small ${\it 125}$}
\psfrag{tnd4a}[]{\small ${\it 97}$}
\psfrag{tnd5a}[]{\small ${\it 79}$}
\psfrag{tnd6a}[]{\small ${\it 67}$}
\psfrag{tnd7a}[]{\small ${\it 58}$}
\psfrag{tnd8a}[]{\small ${\it 51}$}
\psfrag{tnd9a}[]{\small ${\it 46}$}

\psfrag{tne1}[]{\small $3011$}
\psfrag{tne2}[]{\small $1762$}
\psfrag{tne3}[]{\small $1250$}
\psfrag{tne4}[]{\small $968$}
\psfrag{tne5}[]{\small $792$}
\psfrag{tne6}[]{\small $668$}
\psfrag{tne7}[]{\small $580$}
\psfrag{tne8}[]{\small $513$}
\psfrag{tne9}[]{\small $456$}

\psfrag{tne1a}[]{\small ${\it 3010}$}
\psfrag{tne2a}[]{\small ${\it 1761}$}
\psfrag{tne3a}[]{\small ${\it 1249}$}
\psfrag{tne4a}[]{\small ${\it 969}$}
\psfrag{tne5a}[]{\small ${\it 792}$}
\psfrag{tne6a}[]{\small ${\it 669}$}
\psfrag{tne7a}[]{\small ${\it 580}$}
\psfrag{tne8a}[]{\small ${\it 512}$}
\psfrag{tne9a}[]{\small ${\it 458}$}

\psfrag{tcc1}[]{$N=10^2$}
\psfrag{tdd1}[]{$N=10^3$}
\psfrag{tee1}[]{$N=10^4$}
\psfrag{tgg1}[]{\small  exact BL} 
\begin{center}
\includegraphics{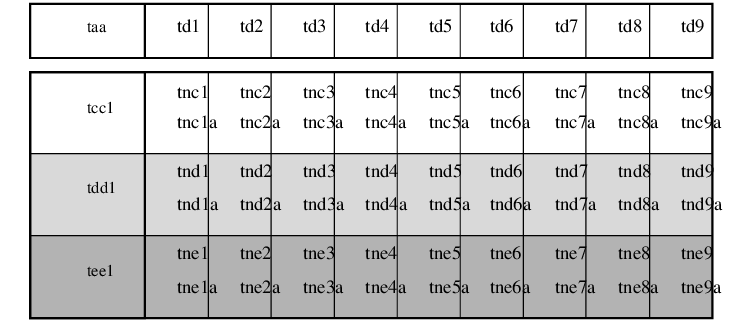}
%
%  magnification 1.00
%
\caption{Occurrences of $D_1$ for the Fibonacci numbers $F_1, F_2,
  \ldots , F_N$, for three different $N$. For comparison, the corresponding best approximation
  of the ``Benford vector'' $N \bigl(\log 2,  \log \frac32 , \ldots , \log
  \frac{10}{9}\bigr)\in \R^9$ by 
  $(N_1, N_2, \ldots, N_9)\in \N^9$ with $\sum_{j=1}^9 N_j=N$ is shown in {\em italics}.}\label{fig1a}
\end{center}
\end{figure}

\begin{example} \label{t14e1a}
For the Fibonacci recursion $x_n = x_{n-1} + x_{n-2}$, $n\ge 3$,
the root of $z^2 = z + 1$ with the largest real part is $\zeta =
\frac12 (1+\sqrt{5})$. As in Example \ref{t14e1}, therefore, $(x_n)$ is Benford unless
$x_1 = x_2 = 0$. In particular, the sequence $(F_n)$ of Fibonacci
numbers is Benford, and so is the sequence of Lucas numbers. As has
been reported repeatedly, the distribution of $D_1$ for $F_1, F_2,
\ldots, F_N$ matches the right-hand side of (\ref{eqNxN1}) remarkably
well even for not-so-large $N$; see, e.g., \cite{BR24, Cal} and Figure
\ref{fig1a}.
\end{example}

A prominent class of matrices amenable to Theorem \ref{thm14} are
irreducible and aperiodic stochastic matrices, often encountered in
homogeneous finite-state Markov chains. As the following result shows,
conformance to BL is prevalent in this context.

\begin{theorem}\label{thm14c}
  Let $d\in \N \setminus \{1\}$, and assume that the random stochastic
  $d\times d$-matrix $P$ has an absolutely continuous distribution
  (w.r.t.\ Lebesgue measure). Then, with probability one, $P$ is
  irreducible and aperiodic. Moreover, with $P^* = \lim_{n\to \infty}
  P^n$ and probability one, the sequences $([P^{n+1} - P^n]_{k \ell})$ and
  $([P^n - P^*]_{k \ell})$ both are Benford or eventually zero for every
  $k,\ell \in \{1,2, \ldots , d\}$.
\end{theorem}

\begin{proof}
See \cite[Thm.\ 12]{BHKR}.
\end{proof}

\begin{example}\label{ex14d}
Let $X,Y$ be independent $U(0,1)$-random variables. By Theorem
\ref{thm14c}, for the two-state Markov chain with one-step transition
matrix
$$
P = \left[
  \begin{array}{cc}
   1- X & X \\ Y & 1- Y
  \end{array}
    \right]\in \R^{2\times 2} \, ,
$$
with probability one, $P$ is irreducible and aperiodic,
$$
P^* = \lim\nolimits_{n\to \infty} P^n = \frac1{X+Y} \left[
  \begin{array}{cc}
  Y & X \\ Y & X
  \end{array}
    \right] \, ,
$$
and the sequences $([P^{n+1} - P^n]_{k \ell})$ and $([P^n - P^*]_{k \ell})$
both are Benford for every $k,\ell \in \{1,2\}$.    
\end{example}
Without the assumption, in Theorem \ref{thm14}, that $A\in \bbs$
satisfy $A^m>0$ for some $m\in \N$, the
Benford properties of the sequences $(|A^n|)$ and $([A^n]_{k \ell})$ are
considerably more delicate. For a comprehensive analysis of the
latter, as well as for studies on BL in more general multi-dimensional
dynamical systems, the reader is referred, e.g., to \cite{BE, BerEsh}.

For a second scenario, let $\bbs$ be the set of all (real-valued)
random variables (on a fixed probability space). Clearly, if $X,Y\in
\bbs$ then also $XY\in \bbs$, and hence in particular $X^n\in \bbs$
for every $n\in \N$.

\begin{theorem}\label{thm14e}
Let $X$ be a (real-valued) random variable. Then:
\begin{enumerate}
\item $(X^n)$ converges in distribution to BL if and only if
  $ \mathbb{E} e^{2\pi i n \log|X|} \to 0$ as $n\to \infty$;
\item $(X^n)$ is Benford with probability one if and only if $\log
  |X|$ is irrational with probability one.
\end{enumerate}
\end{theorem}

\begin{proof}
See \cite[Thm.\ 8.5 and 8.6]{BerAH15}.
\end{proof}

\begin{corollary}\label{cor14f}
If the random variable $X$ has a density (w.r.t.\ Lebesgue measure)
then $(X^n)$ converges in distribution to BL and is Benford with
probability one.
\end{corollary}

\begin{example}\label{ex14g}
If the random variable $X$ is uniformly distributed (on some
non-degenerate interval $\I\subset \R$), exponentially distributed, or
normally distributed (with variance $\sigma^2>0$) then the sequence
$(X,X^2, X^3, \ldots)$ converges in distribution to BL and is Benford
with probability one. By contrast, if $X=10^Y$ and $Y$ is uniformly
distributed on the classical Cantor middle thirds set, then $(X^n)$
does not converge in distribution to BL, even though $X$ is
continuous, i.e., the event $\{X = x\}$ has probability zero for every
$x\in \R$.
\end{example}

For a third and final scenario, let $X$ again be a (real-valued)
random variable, and consider the iid sequence $(X_n)$ with $X_1 =
X$. Clearly, if $X,Y$ are independent then $(X_n Y)$ is an iid
sequence as well. Say that $X$ is {\bf non-discrete} if the event
$\{X\in C\}$ has probability $<1$ for every countable set $C\subset
\R$.

\begin{theorem}\label{thm14h}
Let $X_1, X_2, X_3, \ldots$ be iid random variables. If $X_1$ is
non-discrete then the sequence $(X_1, X_1X_2, X_1X_2X_3, \ldots)$
converges in distribution to BL and is Benford with probability one.
\end{theorem}

\begin{proof}
See \cite[Cor.\ 8.19]{BerAH15}; to appreciate how the conclusion
ultimately hinges on an ergodic theorem, see
also \cite{BerEv} and \cite[Thm.\ 4.24]{BerAH15}.
\end{proof}

For a continuous-time variant of Theorem \ref{thm14h}, recall that a
{\em Brownian motion\/} (or {\em Wiener process\/}) is a continuous stochastic process
on $\R$ with independent and stationary increments that are normally
distributed. Widely used as a model, e.g., for stock prices, a
{\em geometric Brownian motion\/} (or {\em Black--Scholes process\/}) is a stochastic
dynamical system on $\R_+$ whose logarithm is a Brownian motion \cite{EelE76}.

\begin{theorem}\label{thm14i}
Let $(X_t)_{t\ge 0}$ be a geometric Brownian motion. Then $(X_t)_{t\ge
0}$ converges weakly to BL and is Benford with
probability one.
\end{theorem}

\begin{proof}
See \cite[Ex.\ 3.3]{BerEv}, \cite[Thms.\ 2.4 and 3.11]{Sch08} as well
as \cite[Cor.\ 8.19]{BerAH15}.
\end{proof}

%%%%%%%%%%%%%%%%%%%%%%%%%%%%%%%%%%%%%%%%%%%%%%%%%%%%%%%%%%%%%%%%%%%%%%%%%%%%%%%%

\section{A new theorem}\label{sec6}

This final section presents one new result that may help to further
illustrate the ubiquity of BL in dynamical systems. As in all other
sections of this survey, the result is presented in a simple
concrete form rather than in the greatest possible generality. Unlike
in other sections, however, a detailed proof is given here. This is not
only because no proof is available yet in the literature but also, and
perhaps more importantly, because from the proof the reader will get a
clear idea of some of the tools naturally employed for the study of BL
in dynamical systems.

To motivate the result, recall from Section \ref{sec5} that one of the simplest, most
prominent examples of a Benford sequence is any non-trivial solution of the linear
two-step recursion
\begin{equation}\label{eq51}
x_n = x_{n-1} + x_{n-2} \qquad \forall n \ge 3 \, .
\end{equation}
With a view towards applications, e.g., in economics \cite{BD}, it is
natural to consider non-linear variants of (\ref{eq51}) such as, for instance,
\begin{equation}\label{eq52}
x_n = a_1  x_{n-1}^{b_1} + a_2 x_{n-2}^{b_2} \qquad \forall n \ge 3 \, ,
\end{equation}
with parameters $a_1, a_2 \in \R_+$ and $b_1, b_2  >1$. The
dynamics on $\bbs = \R^2_+:= \R_+\times \R_+$ of recurrences like
(\ref{eq52}) have been
studied systematically in \cite{BD}, and will be recalled below to the
extent necessary for the present article. To concisely state the main
result, and in analogy to previous sections, let 
$$
\B  = \bigl\{(x_1, x_2 ) \in \X :  (x_n) \: \mbox{\rm is Benford} \bigr\} \, ,
$$
where $(x_n)$ denotes the unique solution of (\ref{eq52}), given
$(x_1,x_2)$. As it turns out, this solution is Benford for almost all, but
not all, choices of $(x_1, x_2)$.

\begin{theorem}\label{thm19a}
  Let $a_1, a_2 \in \R_+$ and $b_1, b_2 >1$. Then $\B$ has full
  measure but is meagre in $\X$.
\end{theorem}

Though straightforward to state, Theorem \ref{thm19a} requires
for its proof a few preparations, aimed primarily at connecting it to
dynamical systems in a meaningful way. To this end, consider the map
$T:\X \to \X$ given by
$$
T(u,v) = (v, a_2 u^{b_2} + a_1 v^{b_1}) \qquad \forall (u,v) \in
\X \, .
$$
In fact, $T$ is a diffeomorphism from $\X$ onto $T(\X)\subset \X$, and
clearly
$$
(x_{n-1}, x_n) = T(x_{n-2}, x_{n-1}) = T^{n-2} (x_1, x_2) \qquad
\forall n \ge 3 \, .
$$
While it is common practice to thus associate a two-dimensional map
with a two-step recursion like (\ref{eq52}), the reader may want to
appreciate that such an approach does not per se guarantee a
satisfactory understanding of the latter. This merely reflects the
well-documented fact that non-linear two- or multi-step recursions may exhibit an
immense variety of dynamical behaviors for which no
comprehensive theory exists at present; see, e.g., \cite{BD} and the extensive
literature cited there. What makes an approach to Theorem
\ref{thm19a} utilizing $T$ viable, then, is that the dynamics of
(\ref{eq52}) actually is quite simple; see \cite[Thm.\ 2.2 and
Exp.\ 4.2]{BD}. 

\begin{proposition}\label{prop36a}
Let $a_1, a_2\in \R_+$ and $b_1, b_2 > 1$. Then there exists a unique
non-empty, open and bounded set $\A_0\subset \X$ with the following properties:
\begin{enumerate}
\item $\A_0$ is convex, and $\A_0 \supset \{x\in
  \X: x_1^2 + x_2^2 < r^2\}$ for some $r>0$;
  \item the boundary $\partial \A_0$ of $\A_0$ is smooth, and if
    $(x_1, x_2)\in \partial \A_0$ then $(x_n)$ is asymptotically
    $2$-periodic, i.e., $(x_{2n-1})$ and $(x_{2n})$ both converge;
  \item if $(x_1,x_2) \in \A_0$ then $\lim_{n\to \infty}x_n = 0$;
    \item if $(x_1,x_2)\in \X \setminus (\A_0 \cup \partial \A_0)$
      then $\lim_{n\to \infty} x_n = \infty$.
\end{enumerate}
\end{proposition}

Informally put, Proposition \ref{prop36a} asserts that the dynamics
of (\ref{eq52}) neatly partitions the quadrant $\X$ into precisely
three parts: the domains of attraction of $0$ and $\infty$, given by
the non-empty open sets $\A_0$ and $\A_{\infty}:=  \X \setminus (\A_0 \cup
\partial \A_0)$, respectively, as well as their common boundary
$\partial \A_0 = \partial \A_{\infty}$ which is a smooth
curve. Regarding assertion (ii), note that if $(x_1, x_2)\in \partial \A_0$ and
$(x_{2n-1}, x_{2n})\to (p,q) \in \partial \A_0$ then $T(p,q) = (q,p)$
and $T^2(p,q) = (p,q)$, with $p=q$ and $p\ne q$ both being
possibilities, depending on the values of the parameters $a_1,a_2,
b_1, b_2$; see also Figure \ref{fig2} and Example \ref{exPP1} below.

\begin{figure}[ht]
\psfrag{tx1}[]{$x_1$}
  \psfrag{tx2}[]{$x_2$}
\psfrag{tainf}[]{$\A_{\infty}$}
\psfrag{txntoinf}[]{$x_n \to \infty$}
\psfrag{ta0}[]{$\A_0$}
\psfrag{txnto0}[]{$x_n \to 0$}
\psfrag{th12}[]{$\frac12$}
\psfrag{tv12}[]{$\frac12$}
\psfrag{tt2xa}[]{$T^2 (\R^2\! \setminus \! \{(0,0)\})$}
\psfrag{txr2}[]{$T^2 (\R^2_+)$}
\psfrag{txr4}[]{$T^4 (\R^2_+)\subset T^2 (\R^2_+)$}
\psfrag{tt2xb}[]{$= T^2(\R_+^2)$}
\psfrag{tpar}[]{$\partial \A_0 = \partial \A_{\infty}$}
\psfrag{txnto12}[]{$x_n \to \frac12$}
\psfrag{tcap}[]{$(x_{2n-1}, x_{2n}) = T^{2(n-1)}(x_1,x_2)$}
\begin{center}
\includegraphics{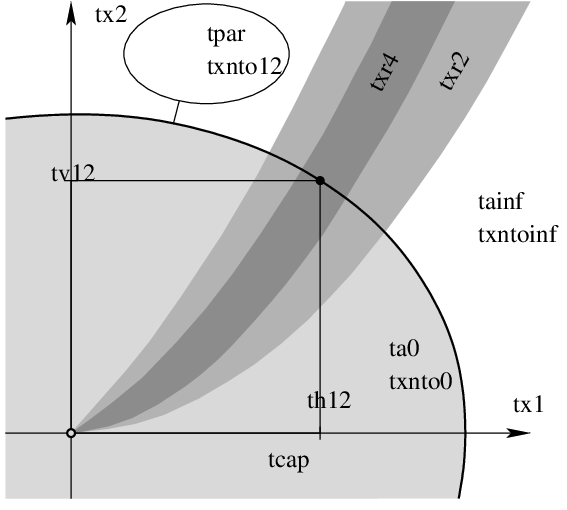}
%
%  magnification 1.00
%
\caption{For the ``quadratic Fibonacci'' recursion of Example
  \ref{exPP1}, that is, for $x_n = x_{n-1}^2 +
  x_{n-2}^2$, $n\ge 3$, the set $\A_0 = \{(x_1, x_2)\in \R^2: \lim_{n\to
    \infty}x_n=0 \}$ is
  non-empty, open, bounded and convex. Its boundary $\partial \A_0$ is
  a smooth oval containing the (unique) fixed point $(\frac12, \frac12)$ of $T$.
If $(x_1,x_2)\in \partial \A_0$ then $\lim_{n\to \infty}x_n =
\frac12$, and clearly $(x_n)$ is {\em not\/} Benford in this case.
By contrast, Corollary \ref{cor19b} shows that $(x_n)$ is Benford for
almost all, but not all, $(x_1,x_2)\in \A_0\cup \A_{\infty}$.
The figure only displays slightly more than one quarter of $\A_0$
(light grey) and $\A_{\infty}$ (white), in light of their obvious
symmetries; note, however, that neither set is symmetric w.r.t.\ the line $x_1=x_2$.}\label{fig2}
\end{center}
\end{figure}

In addition to Proposition \ref{prop36a}, the proof of Theorem
\ref{thm19a} makes use of the following elementary calculus fact
which may also be viewed as a special case of a
{\em shadowing lemma}; see, e.g., \cite[Lem.\ 6.24]{BerAH15}.

\begin{proposition}\label{prop36b}
  Let $(y_n)$ be a sequence in $\R$, and $b>1$. Assume that $\sup_{n\ge
    2}|y_n -  b  y_{n-1}|<\infty$. Then $\sup_{n\in \N}|y_n -
  b^n y |<\infty$ with $y\in \R$ if and only if
$$
y= \widehat{y} := \frac{y_1}{b } + \sum\nolimits_{k=2}^{\infty}
\frac{y_k - b  y_{k-1}}{b^k} \, .
$$
Moreover, if $\lim_{n\to \infty} (y_n -b  y_{n-1})=c$ then 
$\lim_{n\to \infty} (y_n -b^n \widehat{y})=-c/(b-1)$.
\end{proposition}

\begin{proof}[Proof of Theorem \ref{thm19a}]
By Proposition \ref{prop36a}, except for the no\-where dense null\-set
$\partial \A_0$, the two non-empty open sets $\A_0$ and $\A_{\infty}$
make up all of $\X$. Thus it suffices to show that $\B$
has full measure but is meagre in each of these two sets. A detailed argument is presented
here for $(x_1,x_2)\in \A_{\infty}$, that is, for the case of $x_n \to
\infty$. The argument for $(x_1, x_2)\in \A_0$ is completely
analogous and will be commented on only briefly towards the very end of the
proof, with all details being left to the interested reader.

Assume from now on that $(x_1,x_2)\in
\A_{\infty}$. Recall from Theorem \ref{thm1a} that replacing $(x_n)$
by $(a x_n)$ with any $a\in \R_+$ does not affect the Benford
property. When rescaled in this way, (\ref{eq52}) reads
$$
x_n = a_1 a^{b_1-1} x_{n-1}^{b_1} + a_2 a^{b_2 - 1} x_{n-2}^{b_2}
\qquad \forall n \ge 3 \, .
$$
Thus, taking $a= a_1^{-1/(b_1-1)}$ or $a=a_2^{-1/(b_2-1)}$, it can be
assumed that $a_1=1$ or $a_2=1$, respectively. In either case, clearly
$(x_{n-1}, x_n)\in \: ]1,\infty[^2$ for all 
sufficiently large $n\in \N$, and hence $\A_{\infty} = \bigcup_{n\in
  \N} T^{-n} (]1,\infty[^2)$. Since $T$ is a diffeomorphism, it
suffices to show that $\B$ has full measure but is meagre
in $]1,\infty[^2$. This will now be done
via a careful analysis of the sequence $(y_n)$, where $y_n
=\log x_n$ for every $n\in \N$. By Theorem \ref{thm1a}, $(x_n)$ is Benford if and only if $(y_n)$ is
u.d.\ mod 1. Note that $y_n \to \infty$ as well, and 
$$
\U:= \bigl\{ (y_1,y_2)\in \R_+^2 : (y_n) \: \mbox{\rm is 
  u.d.\ mod } 1\bigr\} \subset \R_+^2 
$$
is a Borel set. So, what really will be shown is that $\U$ has full
measure but is meagre in $\R_+^2$. For the convenience of the reader,
three disjoint cases regarding the possible values of the parameters
$b_1, b_2 > 1$ are going to be considered separately. 

\medskip

{\bf Case I:} $b_1^2 > b_2$. As will become clear shortly, in this case the dynamics
of (\ref{eq52}) is dominated by the term $a_1 x_{n-1}^{b_1}$. Assume w.l.o.g.\ that $a_1=1$, i.e.,
consider
\begin{equation}\label{eqNN1}
x_n = x_{n-1}^{b_1} + a_2 x_{n-2}^{b_2} \qquad \forall n \ge 3 \, .
\end{equation}
Note that $x_n>1$ for every $n$, provided that
$x_1,x_2>1$, and then $x_n \to \infty$. In this situation, deduce from $x_n > x_{n-1}^{b_1}$ and
$b_1^2 > b_2$ that
$$
\frac{x_n^{b_1}}{x_{n-1}^{b_2}} > x_{n-1}^{b_1^2 - b_2} \to
\infty \quad \mbox{\rm as } n\to \infty \, .
$$
With the abbreviation
$$
\delta_n  = \delta_n (y_1,y_2):= b_2 y_{n-1} - b_1 y_n \qquad \forall n \ge 2 \, ,
$$
therefore, $\delta_n \le - (b_1^2 - b_2) y_{n-1} <0$ for $n\ge 3$, and
$\delta_n \to - \infty$ whenever
$y_1, y_2 > 0$. Note that $(y_n)$ is given recursively as
$$
y_n =  b_1 y_{n-1} + \log \bigl(
1 + a_2    10^{\delta_{n-1}}
  \bigr) \qquad \forall n \ge 3 \, .
$$
By Proposition \ref{prop36b}, $\lim_{n\to \infty} \bigl(y_n -
b_1^{n-2} h(y ) \bigr)=0$, where for every $y \in
\R_+^2$,
\begin{equation}\label{eqNN2}
h(y) = h(y_1,y_2) = y_2 + \sum\nolimits_{k=1}^{\infty} b_1^{-k} \log
  \bigl(1+ a_2 10^{\delta_{k+1}} \bigr)  \, .
\end{equation}
The function $h:\R_+^2 \to \R$ thus defined is continuous and positive. 
Moreover, $(y_n)$ is u.d.\ mod 1 if and only if $\bigl( b_1^{n-2}
h(y)\bigr)$ is. As will now be shown, the latter is the case precisely
for all $y \in \R_+^2$ belonging to a meagre set of full
measure. Since that set equals $\U$ by definition, this will
establish the desired properties of $\U$.

To this end, note first that $\lim_{n\to \infty} \delta_n b_1^{-n} = - (b_1^2
- b_2) b_1^{-3} h(y)<0$, so $\delta_n \to - \infty$ exponentially
fast as $n\to \infty$. Also,
$$
\delta_n =\delta_n(y) = b_2 y_{n-1} - b_1 y_n = b_2 f^{n-2} (y )_1 -
b_1 f^{n-2} (y )_2 \qquad \forall n \ge 2 \, ,
$$
with the smooth map $f:\R^2 \to \R^2$ given by
\begin{align*}
f  (u,v) & = \bigl(\log T(10^{u},10^{v})_1, \log
         T(10^{u},10^{v})_2 \bigr)
  \\
  & =
\bigl( v ,b_1 v   +
\log ( 1 + a_2  10^{b_2 u - b_1 v} )  \bigr) \, .
\end{align*}
Termwise partial differentiation of
(\ref{eqNN2}) yields, at least formally,
\begin{equation}\label{eqNN3}
\frac{\partial h}{\partial y_2}  = \frac1{1 + a_2   10^{\delta_2}} + \sum\nolimits_{k=1}^{\infty} b_1^{-(k+1)}
\frac{a_2  10^{\delta_{k+2}}}{1+a_2  10^{\delta_{k+2}}} \cdot 
\frac{\partial \delta_{k+2}}{\partial y_2} \, ,
\end{equation}
where, for every $k\ge 1$,
\begin{align*}
\frac{\partial \delta_{k+2}}{\partial y_2} & = b_2
\frac{\partial}{\partial y_2} f^k (y)_1- b_1 \frac{\partial}{\partial y_2} f^k (y)_2 \\
  & =
b_2 \left(
\prod_{\ell =0}^{k-1} Df \circ f^{\ell }(y)
  \right)_{\!\! 21} - b_1 \left(
\prod_{\ell =0}^{k-1} Df \circ f^{\ell }(y)
  \right)_{\!\! 22} \, .
\end{align*}
Together with
$$
Df\circ f^{\ell }(y) = \left[
  \begin{array}{cc} 0 &  \displaystyle  \frac{a_2b_2 10^{\delta_{\ell +2}}}{1 + a_2
    10^{\delta_{\ell +2}}}   \\[4mm]
   1 &  \displaystyle  \frac{b_1}{1 + a_2 10^{\delta_{\ell +2}}} 
        \end{array}
  \right] \qquad \forall \ell  \ge 0 \, ,
  $$
and consequently $|Df \circ f^{\ell }|\le \sqrt{1+b_1^2 + b_2^2}$, this
shows that the series on the right in (\ref{eqNN3}) does converge,
in fact absolutely as well as uniformly in $y\in \R_+^2$. Hence $\partial h/\partial
y_2$ is continuous on $\R_+^2$, and so is $\partial h/\partial y_1$,
by a virtually identical argument. Thus $h$ is $C^1$ on $\R_+^2$.
Moreover, deduce from (\ref{eqNN3}) that by taking $c\in \R_+$
sufficiently large, one can guarantee that $\partial h/\partial
y_2(y)>\frac12$ whenever $\delta_2 < - c$ or, equivalently,
$b_1 y_2 - b_2 y_1 > c$. (Note that $c$ may depend on $a_2 ,b_1,
b_2$, but is independent of $y_1,y_2$; also, the value $\frac12$ for
the lower bound on $\partial h/\partial y_2$ is arbitrary and could be
replaced by any positive number $<1$.) For convenience, let
$$
\V =  \bigl\{
y\in \R_+^2 : b_1 y_2 - b_2 y_1 > c
\bigr\} = \bigl\{
y\in \R_+^2 : \delta_2 < - c
\bigr\}  \, .
$$
Now, recall that $\lim_{n\to \infty}\bigl(y_n - b_1^{n-2} h(y) \bigr)
= 0$, and deduce from Proposition \ref{prop25a} that $(b_1^{n-2} z)$ is u.d.\ mod $1$ for almost all $z\in
\R$, but the set $\bigl\{z\in \R: (b_1^{n-2} z) \:\,  \mbox{\rm is u.d.\
  mod } 1\bigr\}$ is meagre (in $\R$).
Since $\partial h/\partial y_2>\frac12 $ on $\V$, this implies that
for every $y_1\in \R_+$ the $y_1$-sections of $\V \setminus \U$ and $\V\cap\U$, i.e., the sets
\begin{align*}
(\V \setminus \U)_{y_1} & =  \bigl\{ y_2 \in \R_+ : (y_1,y_2)\in \V , \,
                     (y_n) \: \mbox{\rm is {\em not\/} u.d.\ mod } 1\bigr\} \, , \\
  (\V \cap \U)_{y_1} & =  \bigl\{ y_2 \in \R_+ : (y_1,y_2)\in \V , \,
                     (y_n) \: \mbox{\rm is u.d.\ mod } 1\bigr\} \, ,
\end{align*}
are a nullset and meagre in $](c+b_2y_1)/b_1, \infty[$, respectively. This is
because $h(y_1, \cdot)$ is a local diffeomorphism on that open
interval, and hence maps nullsets and meagre sets to nullsets and
meagre sets, respectively. By Fubini's theorem,
$\V\setminus \U$ is a nullset as well, and by the
Kuratowski--Ulam theorem \cite[Thm.\ 15.4]{Ox}, $\V\cap \U$ is meagre in $\V$.
Recall that $\delta_n \to -
\infty$ as $n\to \infty$ for every $y\in \R_+^2$, so $f^n(y)\in \V$ for
all sufficiently large $n$. In other words, $\R_+^2 \subset \bigcup_{n\in \N}
f^{-n} (\V)$. On the one hand, observe that $f(\R_+^2 \setminus \U)
\subset \R_+^2 \setminus \U$, and hence
$$
\R_+^2 \setminus \U  = \bigcup_{n\in \N} \bigl( (\R_+^2 \setminus \U)
\cap f^{-n} (\V) \bigr) \subset \bigcup_{n\in \N} f^{-n} \bigl(
(\R_+^2 \setminus \U)\cap
\V\bigr) =  \bigcup_{n\in \N} f^{-n} (\V \setminus \U)\, .
$$
On the other hand, letting $\V \cap \U \subset\bigcup_{m\in \N} V_m$, where
$V_m \subset \V$ is nowhere dense for each $m\in \N$, observe that
$f(\U)\subset \U$, and hence
$$
\U = \bigcup_{n\in \N} \bigl( \U \cap f^{-n} (\V)\bigr) \subset
\bigcup_{n\in \N} f^{-n} (\V \cap \U) \subset \bigcup_{m,n\in \N}
f^{-n} (V_m)\, .
$$
Since $f$ is a diffeomorphism of $\R_+^2$ onto $f(\R_+^2)$, it follows
that the set
$\U$ has full measure but is meagre in $\R_+^2$. Consequently, 
$\B =  \bigl\{(10^{y_1}, 10^{y_2}) : y\in \U  \bigr\}$
has full measure but is meagre in $]1,\infty[^2$ as well. As
explained earlier, this proves the assertion of the theorem in
case $(x_1,x_2) \in \A_{\infty}$ and $b_1^2>b_2$.

\medskip

{\bf Case II:} $b_1^2 = b_2$. In this case, both terms on the right in (\ref{eq52})
contribute equally to the dynamics. As in Case I, assume w.l.o.g.\ that
$a_1 =1$. Rewriting (\ref{eqNN1}) in terms of the positive quantity $r_n:= x_n /
x_{n-1}^{b_1}$, $n\ge 2$, yields
$$
r_n = \frac{x_n}{x_{n-1}^{b_1}} = 1 + a_2 \left(
\frac{x_{n-2}^{b_1}}{x_{n-1}}
  \right)^{b_1} = 1 + \frac{a_2}{r_{n-1}^{b_1}} \qquad \forall
  n \ge 3 \, .
$$
Thus $r_n = R(r_{n-1}) = R^{n-2} (r_2)$, with the smooth map
$R: \R_+ \to \R_+$ given by $R(r) = 1 + a_2 r^{-b_1}$. Let $\Phi: \X
\to \X$ be the diffeomorphism with
$$
\Phi(u,v) = \bigl( (v/u)^{1/b_1}, v \bigr)
\qquad \forall x \in \X \, .
$$
By Theorem \ref{t8}, with $a_n = R^{n-2}(r)$ and $b_n \equiv b_1$, for
every $r\in \R_+$ the $r$-section of $\Phi^{-1} (\X \setminus \B)$,
i.e., the set
\begin{align*}
\bigl( \Phi^{-1}(\X \setminus \B)\bigr)_{r} & = \bigl\{ x_2 \in \R_+ : \Phi (r, x_2)
                                 \in \X \setminus \B\bigr\} \\
  & = \bigl\{ x_2 \in \R_+ :  x_1 =
(x_2/r)^{1/b_1}, \, (x_n) \: \mbox{\rm is {\em not\/} Benford} \bigr\} \, ,
\end{align*}
is a nullset, whereas
$$
\Phi^{-1} (\B)_r = \bigl\{ x_2 \in \R_+ :  x_1 =
(x_2/r)^{1/b_1}, \, (x_n) \: \mbox{\rm is Benford} \bigr\} 
$$
is meagre (in $\R_+$). Thus
$\Phi^{-1} (\X \setminus \B)$ is a nullset, whereas $\Phi^{-1}(\B)$ is
meagre. Since $\Phi$ is a diffeomorphism, $\B$ has full
measure but is meagre in $\X$.

Note that unlike in Case I above and Case III below, this argument does not require separate
analyses for $\A_{\infty}$ and $\A_0$, nor does it involve taking
logarithms. There is a dynamical reason for this relative simplicity: Since
$$
\Phi^{-1} \circ T \circ \Phi (u,v) = (R(u), R(u) v^{b_1}) \qquad
\forall (u,v) \in \X \, ,
$$
the two-dimensional map $T$ is, up to the diffeomorphism $\Phi$, a {\em skew-product\/}
over the (decreasing) one-dimensional map $R$.

\medskip

{\bf Case III:} $b_1^2 < b_2$. In this case, the dynamics of
(\ref{eq52}) is dominated by the term $a_2 x_{n-2}^{b_2}$. Assume
w.l.o.g.\ that $a_2 =1$, i.e., consider
\begin{equation}\label{eqNN4}
x_n = a_1 x_{n-1}^{b_1} + x_{n-2}^{b_2} \qquad \forall n \ge 3 \, .
\end{equation}
Using the same abbreviation $\delta_n=\delta_n(y_1,y_2) =
b_2 y_{n-1} - b_1 y_n$ for $n\ge 2$ as in Case I, note that
\begin{equation}\label{eqNN5}
  y_n = b_2 y_{n-2} + \log \bigl(
1 + a_1    10^{- \delta_{n-1}}
  \bigr) \qquad \forall n \ge 3 \, . 
\end{equation}
With the goal of applying Proposition \ref{prop36b} in a similar
way as in Case I, it will first be shown that, given any $y\in \R_+^2$,
\begin{equation}\label{eqNN6}
(\delta_{2n-1}), \, (\delta_{2n}) \: \mbox{\rm both converge in } \R
\cup \{\infty\} \, .
\end{equation}
As a first step towards establishing (\ref{eqNN6}), deduce from
(\ref{eqNN5}) that
\begin{align}\label{eqNN7}
 \delta_{2n+1} & =  b_2  \log \bigl(
 a_1 +  10^{\delta_{2n-1}}
\bigr) \: - \nonumber \\
& \quad - b_1  \log \bigl(
1 + a_1    10^{- (b_2 - b_1^2)y_{2n-1}} (a_1 + 10^{\delta_{2n-1}})^{b_1}
  \bigr) 
\qquad \forall n \ge 2 \, .
\end{align}
Motivated by (\ref{eqNN7}), for every $\varepsilon\ge 0$ consider the
smooth map $R_{\varepsilon} :\R\to \R$ given by
$$
R_{\varepsilon} (r) = b_2 \log (a_1 + 10^r) - b_1 \log (1 +
\varepsilon (a_1 + 10^r)^{b_1}) \, .
$$
It is readily checked that $0<R_{\varepsilon}'(r) < b_2$ for all $r\in
\R$, and
$$
R_{\varepsilon} (-\infty)  = b_2 \log a_1 - b_1 \log (1 + \varepsilon
a_1^{b_1} ) \, , \quad
R_{\varepsilon}' (-\infty)  = 0 \, .
$$
Since $R_{\varepsilon}$ is strictly increasing, the sequence $\bigl(
R_{\varepsilon}^n (r) \bigr)$ is monotone for every $r\in \R$, and hence
either converges to a fixed point of $R_{\varepsilon}$, or else
$\lim_{n\to \infty} R_{\varepsilon}^n (r)=\infty$. (Note that
$\lim_{n\to \infty} R_{\varepsilon}^n (r) = -\infty$ is impossible
because $R_{\varepsilon}(-\infty)>-\infty$.) Moreover, the map
$R_0$ is strictly convex, with $R'_0(\infty) = b_2 > 1$, and so has
at most two fixed points. Given any $c\in \R_+$, note that
$R_{\varepsilon}\le R_0$ and $R_{\varepsilon}\uparrow R_0$ uniformly
on $]-\infty, c]$ as $\varepsilon \downarrow 0$. Thus,
if $R_0$ does have fixed points and if $c$ is larger than each of these fixed
points, then $R_{\varepsilon}$ has fixed points in $]-\infty, c]$
also, and with $r_{\varepsilon}^-,r^+_{\varepsilon}$ denoting the
smallest and (not necessarily different) largest fixed point of
$R_{\varepsilon}$ in $]-\infty, c]$, respectively, the limits
$\lim_{\varepsilon\downarrow 0} r^-_{\varepsilon}$,
$\lim_{\varepsilon\downarrow 0} r^+_{\varepsilon}$ both exist and are
fixed points of $R_0$. Now, recall that $b_1^2 < b_2$ and $y_n \to
\infty$. Given $\varepsilon>0$, therefore, there exists $N_{\varepsilon}\ge
2$ so that $0< a_1 10^{-(b_2 - b_1^2) y_{2n-3}} < \varepsilon$ for
all $n\ge N_{\varepsilon}$, and (\ref{eqNN7}) yields
$$
R_{\varepsilon} (\delta_{2n-1}) < \delta_{2n+1} < R_0 (\delta_{2n-1})
\qquad \forall n \ge N_{\varepsilon} \, .
$$
It follows that
\begin{equation}\label{eqNN8}
R_{\varepsilon}^n (\delta_{2N_{\varepsilon}-1}) < \delta_{2(n+N_{\varepsilon})-1} < R_0^n (\delta_{2N_{\varepsilon}-1})
\qquad \forall n\in \N \, .
\end{equation}
At this point, it is convenient to consider the three possibilities
regarding the fixed points of $R_0$ individually.

First, if $R_0$ has {\em no fixed point\/} then, given any $c\in \R_+$, take
$\varepsilon>0$ so small that $R_{\varepsilon}$ has no fixed point in
$]-\infty, c]$ either. In this case, letting $n\to \infty$ in (\ref{eqNN8})
yields
$$
c\le \lim\nolimits_{n\to \infty} R_{\varepsilon}^n
(\delta_{2N_{\varepsilon} - 1}) \le \liminf\nolimits_{n\to \infty}
\delta_{2n-1} \, .
$$
Since $c\in \R_+$ has been arbitrary, $\delta_{2n-1} \to \infty$.

Second, if $R_0$ has {\em a single fixed point\/} $r_0$ (which necessarily is
attracting at the left and repelling at the right) then, given any
$c\in \R_+$ with $c>r_0$, take
$\varepsilon>0$ so small that $r_0 < r_{\varepsilon}^+ < c$. From
(\ref{eqNN8}), it is clear that
$$
r_{\varepsilon}^- \le \lim\nolimits_{n\to \infty} R_{\varepsilon}^n
(\delta_{2N_{\varepsilon} -1}) \le \liminf\nolimits_{n\to \infty}
\delta_{2n-1} \, .
$$
If $\limsup_{n\to \infty} \delta_{2n-1}>r_{\varepsilon}^+$ then there
exists $N\ge N_{\varepsilon}$ so that $\delta_{2N-1} >
r_{\varepsilon}^+$, and hence $\delta_{2n-1} > R_{\varepsilon}^{n-N}
(\delta_{2N-1})$ for all $n\ge N$, which in turn yields $\liminf_{n\to
\infty} \delta_{2n-1} \ge c$, and since $c>r_0$ has been arbitrary,
$\delta_{2n-1} \to \infty$. By contrast, if $\limsup_{n\to \infty}
\delta_{2n-1}\le r_{\varepsilon}^+$ then clearly
$$
r_{\varepsilon}^-  \le \liminf\nolimits_{n\to \infty}
\delta_{2n-1} \le \limsup\nolimits_{n\to \infty}
\delta_{2n-1} \le r_{\varepsilon}^+ \, ,
$$
and letting $\varepsilon \downarrow 0$ yields $\lim_{n\to \infty}
\delta_{2n-1} = r_0$.

Third, if $R_0$ has {\em two fixed points\/} $r_0^- < r_0^+$ (where $r_0^-$ and $r_0^+$
necessarily are attracting and repelling, respectively) then pick any
$c> r_0^+$ and take $\varepsilon > 0$ so small that $r_0^+ <
r_{\varepsilon}^+ < c$. On the one hand, if
$\liminf_{n\to \infty} \delta_{2n-1} < r_0^+$ then there exists $N\ge
N_{\varepsilon}$ so that $\delta_{2N-1} < r_0^+$, and hence
$\delta_{2n-1} < R_0^{n-N} (\delta_{2N-1})$ for all $n\ge N$, so
$\limsup_{n\to \infty} \delta_{2n-1} \le \lim_{n\to \infty} R_0^{n-N}
(\delta_{2N-1}) = r_0^-$. Since clearly $\liminf_{n\to \infty}
\delta_{2n-1}\ge r_{\varepsilon}^-$, letting $\varepsilon \downarrow
0$ yields $\lim_{n\to \infty} \delta_{2n-1} = r_0^-$. On the other hand, if
$\limsup_{n\to \infty} \delta_{2n-1} > r_{\varepsilon}^+$ then
$\delta_{2n-1} \to \infty$, just as in the case of a single fixed point. The only remaining
possibility is that
$$
r_0^+ \le  \liminf\nolimits_{n\to \infty}
\delta_{2n-1}  \le  \limsup\nolimits_{n\to \infty}
\delta_{2n-1} \le r_{\varepsilon}^+ \, ,
$$
where letting $\varepsilon \downarrow
0$ yields $\lim_{n\to \infty} \delta_{2n-1} = r_0^+$.

In summary, it has thus been shown that $(\delta_{2n-1})$ converges in
$\R\cup \{\infty\}$ for every choice of $y \in \R_+^2$. In fact, every
{\em finite\/} limit of $(\delta_{2n-1})$ has been identified as a fixed point
of $R_0$. Moreover, it is clear that the set
$$
\V_o:= \bigl\{ y \in \R_+^2 : \lim\nolimits_{n\to \infty}
\delta_{2n-1} = \infty  \bigr\}
$$
is non-empty and open. If $R_0$ has no fixed point then clearly $\V_o= \R_+^2$,
but otherwise $\V_o$ may be smaller.

To establish (\ref{eqNN6}), it remains to consider the sequence
$(\delta_{2n})$. For this, simply notice that (\ref{eqNN7}) remains
valid with $n$ replaced by $n+\frac12$. Thus, the same argument
as above shows that $(\delta_{2n})$ either converges to a finite
limit (which again must be a fixed point of $R_0$) or else
$\delta_{2n} \to \infty$. This fully proves (\ref{eqNN6}). Again, the
set
$$
\V_e:= \bigl\{ y \in \R_+^2 : \lim\nolimits_{n\to \infty}
\delta_{2n} = \infty  \bigr\}
$$
is non-empty and open. In addition, the identity
$$
b_2\delta_{2n-1} + b_1\delta_{2n} = b_2 (b_2 - b_1^2) y_{2n-2} - b_1^2
\log (1 + a_1 10^{-\delta_{2n-1}}) \qquad \forall n \ge 2 \, ,
$$
shows that $\lim_{n\to \infty}\delta_{2n-1}$ and $\lim_{n\to
  \infty}\delta_{2n}$ cannot both be finite. In other words, $\V_o
\cup \V_e = \R_+^2$. 

The scene is now set for a (``shadowing'') argument
utilizing Proposition \ref{prop36b} which is similar to the one in Case
I, yet slightly more involved. Thus, assume for instance that $y\in \V_o$, i.e., $\lim_{n\to \infty}
\delta_{2n-1} = \infty$. As in Case I, deduce from (\ref{eqNN5}) that 
 $\lim_{n\to \infty} \bigl(y_{2n} - b_2^{n-1} h(y) \bigr)=0$, where
 for every $y\in \R_+^2$,
 \begin{equation}\label{eqNN10}
h(y) = y_2 +\sum\nolimits_{k=1}^{\infty} b_2^{-k}  \log \bigl(1 + a_1
  10^{-\delta_{2k+1}}\bigr)  \, ,
\end{equation}
and hence $(y_{2n})$ is u.d.\ mod $1$ if and only if $\bigl( b_2^{n-1}
h(y)\bigr)$ is. In analogy to Case I, it is readily seen that $\lim_{n\to
  \infty} \delta_{2n-1} b_2^{-n}$ is finite and positive, so
$\delta_{2n-1} \to \infty$ exponentially fast as $n\to \infty$. Letting
$$
  \E_o = \left\{
y\in \V_o : (y_{2n}) \:\mbox{\rm or} \left( \frac{b_1}{b_2} \,
  y_{2n}\right)  \:\mbox{\rm is {\em not\/} u.d.\ mod } 1
    \right\} \, ,
$$
it will now be shown that
\begin{equation}\label{eqNP12a}
\E_o \subset \V_o \: \mbox{\rm is a nullset}\, , 
\end{equation}
whereas later it will also be shown that
\begin{equation}\label{eqNP12b}
\V_o \cap \U \: \mbox{\rm is meagre (in $\V_o$)}\, . 
\end{equation}
As a first step towards (\ref{eqNP12a}), and analogously to Case I, observe that
$$
\delta_{2n+1} = b_2 y_{2n} - b_1 y_{2n+1} = b_2 g^{2n-1} (y)_1 - b_1
g^{2n-1} (y)_2 \qquad \forall n \ge 1 \, ,
$$
with the smooth map $g:\R^2 \to \R^2$ given by
$$
g  (u,v) = 
\bigl(v ,b_2 u  +
\log ( 1 + a_1  10^{b_1 v - b_2 u} )  \bigr) \, .
$$
Note that $h:\R_+^2 \to \R$ is continuous, and $h>0$ on
$\V_o$. (Formal) term\-wise partial differentiation of (\ref{eqNN10}) on $\V_o$ yields
\begin{equation}\label{eqNN12}
\frac{\partial h}{\partial y_2}  = 1 -  \sum\nolimits_{k=1}^{\infty} b_2^{-k}
\frac{a_1  10^{-\delta_{2k+1}}}{1+a_1  10^{-\delta_{2k+1}}} \cdot 
\frac{\partial \delta_{2k+1}}{\partial y_2} \, ,
\end{equation}
where, for every $k\ge 1$,
\begin{align*}
\frac{\partial \delta_{2k+1}}{\partial y_2} & = b_2
\frac{\partial}{\partial y_2} g^{2k-1} (y)_1- b_1 \frac{\partial}{\partial y_2} g^{2k-1}(y)_2 \\
  & =
b_2 \left(
\prod_{\ell =0}^{2k-1} Dg \circ  g^{\ell }(y)
  \right)_{\!\! 21} - b_1 \left(
\prod_{\ell =0}^{2k-1} Dg \circ g^{\ell }(y)
  \right)_{\!\! 22} \, .
\end{align*}
Observe now that
$$
Dg\circ g^{\ell }(y) = \left[
  \begin{array}{cc} 0 &  \displaystyle  \frac{b_2}{1 + a_1
    10^{-\delta_{\ell+2}}}   \\[4mm]
   1 &  \displaystyle  \frac{a_1 b_1 10^{-\delta_{\ell +2}}}{1 + a_1
       10^{-\delta_{\ell +2}}} 
        \end{array}
  \right] \qquad \forall \ell  \ge 0 \, .
  $$
Recall that $\delta_{2n-1} \to \infty$, while the sequence
$(\delta_{2n})$ either converges to a finite limit, or else
$\delta_{2n} \to \infty$ as well. In either case, and similarly to Case I,
$\sup_{\ell \ge 0} |Dg \circ g^{\ell}|<\infty$, and since $\delta_{2n-1} \to
\infty$ exponentially fast, (\ref{eqNN12}) shows that $\partial h/\partial
y_2$ is continuous on $\V_o$, as is $\partial h/\partial y_1$. Again,
pick $c\in \R_+$, possibly depending on $a_1,b_1,b_2$ but
independent of $y_1,y_2$, so large that $\partial h/\partial y_2 >
\frac12 $ whenever $\delta_3 > c$. With this, let $\W_o = \{y\in \V_o: \partial h/\partial y_2 >
\frac12 \}$. Since $h(y_1, \cdot)$ is a local diffeomorphism on every
non-empty $y_1$-section of $\W_0$,
the $y_1$-section
$$
(\E_o \cap \W_o)_{y_1} = \left\{y_2 \in \R_+: 
y\in \W_o, \,  (y_{2n}) \:\mbox{\rm or} \left( \frac{b_1}{b_2} \,
  y_{2n}\right)  \:\mbox{\rm is {\em not\/} u.d.\ mod } 1
    \right\} 
$$
is a nullset. Thus, $\E_o \cap \W_o$ is a nullset as well. Since
$\delta_{2n-1} \to \infty$ as $n\to \infty$, $g^{2n}(y)\in \W_o$ for
all sufficiently large $n$. Thus $\V_o \subset \bigcup_{n\in \N}
g^{-2n} (\W_o)$, and since clearly $g^2 (\E_o) \subset \E_o$, 
$$
\E_o = \E_o \cap \bigcup_{n\in \N}g^{-2n} (\W_o) = \bigcup_{n\in \N} \bigl(
\E_o \cap g^{-2n} (\W_o) \bigr) \subset \bigcup_{n\in \N} g^{-2n}
(\E_o \cap \W_o) \, .
$$
Since $g$ is a diffeomorphism of $\R_2^+$ onto $g(\R_+^2)$,
this proves (\ref{eqNP12a}).

Now, to establish (\ref{eqNP12b}), observe first that, as in the proof
of \cite[Thm.\ 6.46]{BerAH15}, given any $y\in \W_o$, $\varepsilon
>0$, and $0<t<\frac14$ so small that
$$
10 t \left(
\frac{\log 2 - \log t}{\log b_2} + 1
\right)
< 1 \, ,
$$
there exists $z\in \W_o$ with $|z-y|<\varepsilon$ and
$$
\liminf\nolimits_{N\to \infty} \frac{\# \{ 1\le n \le N : z_{2n} \in
  [\frac12 - t , \frac12+t] +\Z \} }{N} > 10 t \, ,
$$
where $z_n = \log x_n$ for all $n$, and $(x_n)$ is the solution of
(\ref{eqNN4}) with $(x_1,x_2)=(10^{z_1}, 10^{z_2})$. It follows that
\begin{equation}\label{eqXQ1}
\liminf\nolimits_{N\to \infty} \frac{\# \{ 1\le n \le N : z_{n} \in
  [\frac12 - t , \frac12+t] +\Z \} }{N} > 5 t \, .
\end{equation}
Consider the continuous, $1$-periodic function $\psi:\R \to [0,1]$
with
$$
\psi(s) =\left\{
  \begin{array}{ll}
    0 & \mbox{\rm if }\mbox{\rm dist}(s, \frac12 + \Z) \ge 2t \, , \\[1mm]
    2 - \mbox{\rm dist}(s, \frac12 + \Z) /t & \mbox{\rm if }t  \le
                                              \mbox{\rm dist}(s,
                                              \frac12 + \Z) < 2t \, , \\[1mm]
    1 & \mbox{\rm if }\mbox{\rm dist}(s, \frac12 + \Z) < t \, , 
    \end{array}
  \right.
$$
and note that $\int_0^1 \psi(s)\, {\rm d}s = 3t$. By (\ref{eqXQ1}),
$$
\liminf\nolimits_{N\to \infty} \frac1{N} \sum\nolimits_{n=1}^N \psi
(z_n) > 5t \, ,
$$
and hence
\begin{equation}\label{eqXQ2}
\liminf\nolimits_{N\to \infty} \left| \frac1{N} \sum\nolimits_{n=1}^N \psi
(z_n) - \int_0^1 \psi(s) \, {\rm d}s
\right| > 2 t \, .
\end{equation}
Now, for every $m\in \N$ define $V_m \subset \W_o$ as
$$
V_m = \left\{
y\in \W_o : \left| \frac1{N} \sum\nolimits_{n=1}^N \psi
(y_n) - \int_0^1 \psi(s) \, {\rm d}s
\right| \le t \: \forall N \ge m
  \right\}\, .
$$
Note that the set $V_m$ is closed (in $\W_o$) by continuity, but has empty interior,
by (\ref{eqXQ2}). Also, $V_1\subset V_2 \subset \cdots$, and $\W_o
\cap \U \subset \bigcup_{m\in \N} V_m$ by Weyl's criterion \cite[Cor.\
I.1.2]{KN}. Recall from above that $\V_o \subset \bigcup_{n\in \N}
g^{-2n} (\W_o)$, and clearly $g^2(\U) \subset \U$. Thus
$$
\V_o \cap \U  \subset \bigcup_{n\in \N} (g^{-2n} (\W_o) \cap \U)
               \subset \bigcup_{n\in \N} g^{-2n} (\W_o \cap \U) \subset \bigcup_{m,n\in \N} g^{-2n} (V_m) \, ,
$$
which shows that $\V_o\cap \U$ is meagre, and hence proves (\ref{eqNP12b}).

An entirely analogous argument utilizing the set
$$
\E_e := \left\{
y\in \V_e :  \,  (y_{2n-1}) \:\mbox{\rm or} \left( \frac{b_1}{b_2} \,
  y_{2n-1}\right)  \:\mbox{\rm is {\em not\/} u.d.\ mod } 1
    \right\} \, ,
$$
shows that, similarly to (\ref{eqNP12a}) and (\ref{eqNP12b}),
$$
\E_e\subset \V_e \: \mbox{\rm is a nullset, and $\V_e \cap \U$ is
  meagre}\, .
$$
Clearly, therefore, $\U = (\V_o\cap \U)
\cup (\V_e \cap \U)$ is meagre as well, and it only remains to be shown
that $ \U$ has full measure in $\R_+^2$. To this end, note that $\E_o
\cup \E_e$ is a nullset, and pick any $y\in \R_+^2 \setminus (\E_o \cup
\E_e )$. Assume, for instance, that $y\in \V_o$. Then $(y_{2n})$ and
$(b_1b_2^{-1} y_{2n})$ both are u.d.\ mod $1$. On the one hand, if
$y\in \V_{e}$ then $(y_{2n-1})$ is u.d.\ mod $1$ as well. On the other
hand, if $y\not \in \V_e$ then $(\delta_{2n})$ converges to a finite
limit. In this case, recall that
$$
y_{2n-1} = \frac{b_1}{b_2} y_{2n} + \frac{1}{b_2} \delta_{2n} \qquad
\forall n \ge 1 \, ,
$$
which shows that, again, $(y_{2n-1})$ is u.d.\ mod $1$. In either case,
therefore, $(y_{2n})$ and $(y_{2n-1})$ both are u.d.\ mod $1$, and
hence $(y_n)$ is as well. By completely analogous reasoning, the same
conclusion holds whenever $y\in \V_e$. Thus $(y_n)$ is u.d.\ mod $1$
whenever $y\in \R_+^2 \setminus (\E_o \cup
\E_e )$. In other words, $(\R_+^2 \setminus \U) \subset \E_o \cup
\E_e$, and so $\R_+^2 \setminus \U$ is a nullset. Thus, $\U$ has
full measure but, as seen earlier, is meagre in $\R_+^2$. This
proves the assertion of the theorem in case $(x_1,x_2)\in \A_{\infty}$
and $b_1^2 < b_2$.
 
\medskip

To conclude the overal proof, the case of $(x_1,x_2)\in \A_0$ still needs to
be considered. In essence, the analysis in this
case proceeds exactly as in the case of $(x_1,x_2)\in \A_{\infty}$,
but with Cases I and III interchanging their roles. To see this, fix
$0<\xi < 1$ so that $T (]0,\xi[^2) \subset \: ]0,\xi[^2$. Thus, if
$0<x_1,x_2<\xi$ then $0<x_n < \xi$ for every $n$, and $\lim_{n\to
  \infty} x_n = 0$. In analogy to before, it suffices to show that
$\B$ has full measure but is meagre in $]0,\xi[^2$. Assume for
instance that $b_1^2 < b_2$, and consider (\ref{eq52}) in its
rescaled form (\ref{eqNN1}). With $\delta_n$ defined as before,
$$
10^{\delta_n} = \frac{x_{n-1}^{b_2}}{x_n^{b_1}} < x_{n-1}^{b_2 -
  b_1^2} \to 0 \quad \mbox{\rm as } n \to \infty \, .
$$
Thus, $\delta_n \to -\infty$, and the same technique as in Case I
applies. By contrast, if $b_1^2 > b_2$ then consider 
(\ref{eqNN4}), with (\ref{eqNN7}) remaining valid as before, and again
$(b_2-b_1^2)y_{2n-1} \to \infty$ as $n\to \infty$. Consequently, the
same technique as in Case III applies. Finally, as observed
earlier, the argument in Case II applies equally if $(x_1,x_2)\in \A_0$.
\end{proof}

\begin{corollary}\label{cor19b}
Let $a_1,a_2\in \R_+$. Then the solution $(x_n)$ of the recursion $x_n =
a_1 x_{n-1}^2 + a_2 x_{n-2}^2$, $n\ge 3$, is Benford for almost all
$(x_1,x_2)\in \R^2$, but the set $\{(x_1,x_2)\in \R^2: (x_n) \:
\mbox{\rm is not Benford}\}$ is residual in $\R^2$.
\end{corollary}

\begin{proof}
This follows immediately from Theorem \ref{thm19a}, since $T:\R^2 \to
\R^2$ given by $T(u,v) =  (v, a_2u^2 + a_1v^2)$ is a diffeomorphism when
restricted to any one of the four open quadrants, and $T^2 (\R^2
\setminus \{(0,0)\})\subset \R_+^2$. Notice that Case I in the above proof
of Theorem \ref{thm19a} applies always, since $b_1^2 = 4 > 2 = b_2$.
\end{proof}

\begin{example}\label{exPP1}
  Consider the ``quadratic Fibonacci'' recursion
  \begin{equation}\label{eqPP1}
    x_n = x_{n-1}^2 + x_{n-2}^2 \qquad \forall n \ge 3 \, .
  \end{equation}
  By Corollary \ref{cor19b}, $(x_n)$ is Benford for almost all, but
  not all $(x_1,x_2)\in \R^2$. The boundary of $\A_0$, the domain of
  attraction of $0$, is a smooth oval containing $(\frac12,
  \frac12)$. In fact, $\lim_{n\to \infty} x_n = \frac12$ whenever
  $(x_1, x_2)\in \partial \A_0$; see also Figure \ref{fig2}. Note that $x_n \to \infty$ for every
  $(x_1, x_2)\in \Z^2\setminus \{(0,0)\}$, but the authors do not know
  whether or not any of the resulting
  (doubly-exponentially growing) integer sequences
  are Benford.

  As a slightly modified version of (\ref{eqPP1}), consider the recursion
 $$
    x_n = x_{n-1}^2 + 4 x_{n-2}^2 \qquad \forall n \ge 3 \, .
 $$
  Again, $\partial \A_0$ is a smooth oval, and $(\frac15, \frac15) \in
  \partial \A_0$. Unlike for (\ref{eqPP1}), however, picking $(x_1,
  x_2)\in \partial \A_0$ now yields a divergent sequence $(x_n)$
  unless $(|x_1|, |x_2|) = (\frac15 , \frac15)$ since, for instance,
  $\lim_{n\to \infty} x_{2n-1} = \frac1{30} (5+\sqrt{5})$ and
  $\lim_{n\to \infty} x_{2n} = \frac1{30} (5-\sqrt{5})$ whenever
  $(x_1,  x_2)\in \partial \A_0$ and $|x_1|>\frac15$.
\end{example}

\begin{example}\label{exom}
The conclusion of Theorem \ref{thm19a} may remain valid, at least
partially, even if one of the parameters $b_1, b_2$ is $\le 1$. For
example, it is readily seen that for the recursion
$$
 x_n = x_{n-1}^2 +  \sqrt{x_{n-2}}  \qquad \forall n \ge 3 \, ,
 $$
 $x_n\to \infty$ for every $(x_1 , x_2)\in \X = \R_+^2$, with $\B$ having
 full measure but being meagre in $\X$.

By contrast, for the slightly modified recursion
$$
 x_n = x_{n-1}^2 +  \frac14 \sqrt{x_{n-2}}  \qquad \forall n \ge 3 \, ,
 $$
 it can be shown that $\lim_{n\to \infty} x_n = \zeta^2_-\approx 0.07268$
 whenever $0< x_1,x_2 < \zeta^2_+\approx 0.7015$, where $\zeta_-$ and
 $\zeta_+$ is the smaller and larger positive root of
 $4z^3 - 4z + 1 = 0$, respectively; in particular, $(x_n)$ is {\em
   not\/} Benford for any such $(x_1,x_2)$. Even in this example,
 however, $\A_{\infty} \supset \: ]\zeta_+^2, \infty[^2$, and $\B$ has
 full measure but is meagre in $\A_{\infty}$.
\end{example}

In conclusion, the reader may observe that even though Theorem
\ref{thm19a} does not apply to them, the recursions of Example
\ref{exom} do generate an abundance of Benford sequences. So, indeed, do
many other non-linear recursions and more general dynamical
systems. The study of BL in dynamical systems will undoubtedly remain
an intriguing and multi-faceted area of research for many years to come.  

\section*{Acknowledgements}

The first author was partially supported by an NSERC Discovery
Grant. Both authors are grateful to an anonymous referee for several helpful comments and corrections.

%%%%%%%%%%%%%%%%%%%%%%%%%%%%%%%%%%%%%%%%%%%%%%%%%%%%%%%%%%%%%%%%%%%%%%%%%%%%%%%%

%%%%%%%%%%%%%%%%%%%%%%%%%%%%%%%%%%%%%%%%%%%%%%%%%%%%%%%%%%%%%%%%%%%%%%%%%%%%%%%%

\end{document}